\documentclass{article}
\usepackage{graphicx} 

\usepackage[utf8]{inputenc}
\usepackage{geometry}
\usepackage{amsmath, amssymb}
\usepackage{algorithm}[boxed,boxedruled]
\usepackage[algo2e]{algorithm2e}[linesnumbered,vlined]
\usepackage{mathrsfs}
\usepackage{amsthm}
\usepackage{float}
\usepackage{graphicx}
\usepackage{caption}
\usepackage{subcaption}
\usepackage{authblk}
\usepackage{orcidlink}
\usepackage{xcolor}

\newtheorem{theorem}{Theorem}
\newtheorem*{theorem*}{Theorem}

\begin{document}

\title{A Unifying Family of Data-Adaptive Partitioning Algorithms}

\author[a,\orcidlink{0000-0002-3209-4705}]{Maria Emelianenko \footnote{E-mail: memelian@gmu.edu}}
\author[a,\orcidlink{guyoldaker4@gmail.com}]{Guy B. Oldaker IV\footnote{E-mail: goldaker@gmu.edu} }

\affil[a]{Department of Mathematical Sciences, George Mason University, 4400 University Dr, Fairfax, VA 22030}

\maketitle

\begin{abstract}
Clustering algorithms remain valuable tools for grouping and summarizing the most important aspects of data. Example areas where this is the case include image segmentation, dimension reduction, signals analysis, model order reduction, numerical analysis, and others.  As a consequence, many clustering approaches have been developed to satisfy the unique needs of each particular field.  In this article, we present a family of data-adaptive partitioning algorithms that unifies several well-known methods (e.g., k-means and k-subspaces).  Indexed by a single parameter and employing a common minimization strategy, the algorithms are easy to use and interpret, and scale well to large, high-dimensional problems.  In addition, we develop an adaptive mechanism that (a) exhibits skill at automatically uncovering data structures and problem parameters without any expert knowledge and, (b) can be used to augment other existing methods.  By demonstrating the performance of our methods on examples from disparate fields including subspace clustering, model order reduction, and matrix approximation, we hope to highlight their versatility and  potential for extending the boundaries of existing scientific domains. We believe our family's parametrized structure represents a synergism of algorithms that will foster new developments and directions, not least within the data science community.
\end{abstract}

\section{Introduction}
\begin{figure}[H]
		\begin{center}
        \includegraphics[width=15.4cm, height=6cm]{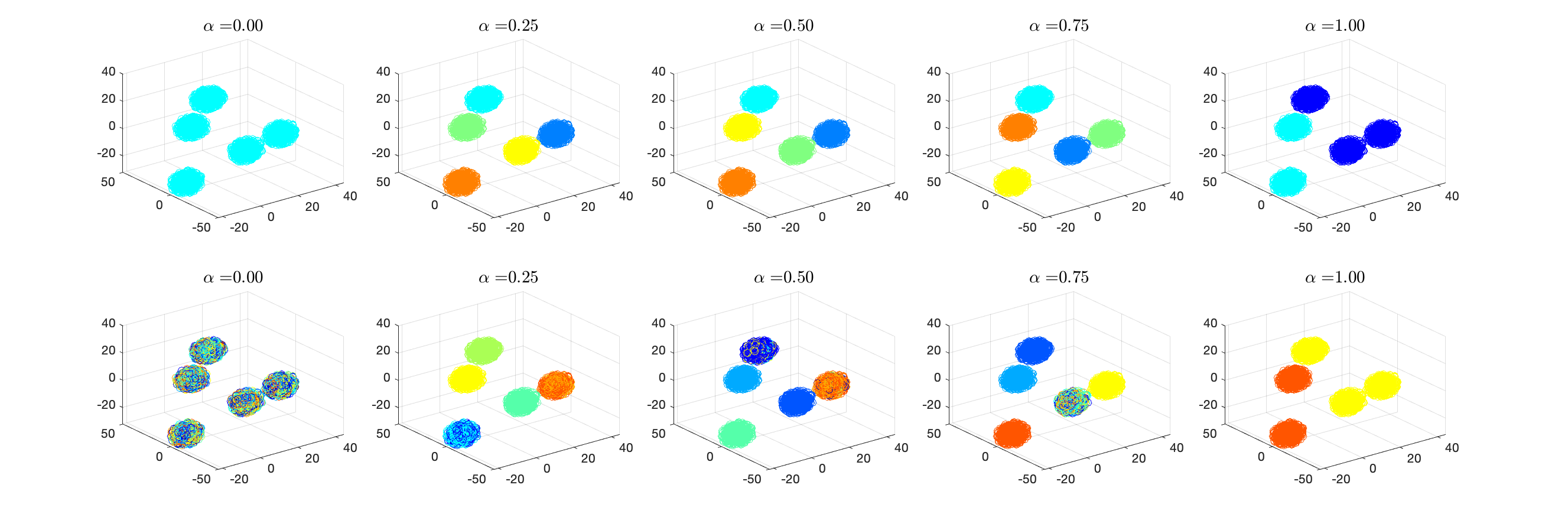}
		\end{center}
		\caption{Example performance of our algorithm on an idealized clustering task with (top) and without (bottom) adaptation as the indexing parameter, $\alpha$, is varied.  The data consist of five multi-variate Gaussian point clouds in $\mathbb{R}^{6000}$ (the data were projected onto $\mathbb{R}^{250}$ via a Gaussian random embedding prior to clustering \cite{dong2021simpler}).  Both algorithms are initialized for $k = 10$ clusters and total dimension 250 (See \ref{section1} for algorithm input details).  Note that for $\alpha \in \{0.25,0.5,0.75\}$ the adaptive variant uncovers the correct number of clusters and labeling, while the non-adaptive version struggles.  The tSNE algorithm \cite{hinton2002stochastic}\cite{van2008visualizing} is used for visualization. }
		\label{fig:spheres}
  \end{figure}

A foundational aspect of data science continues to be the distillation of the essential structure(s) within data.  Examples of this ongoing pursuit can be found among the theories and methods in many diverse application areas including dimension reduction, manifold learning, matrix approximation, signals analysis, and clustering.  In this last area, the objective is to partition a given data set into groups of similar objects.  Applications include computer vision, gene expression analysis, image segmentation, and fluid dynamics \cite{zhang2012hybrid}\cite{kaiser2014cluster}\cite{wang2022convergence}\cite{petegrosso2020machine}.  In addition to organizing data into disjoint groups, clustering provides the practitioner with a means of efficiently summarizing the contents of each partition member.  Techniques/results of this kind have seen extensive use in dimension reduction \cite{kerschen2002non}\cite{kerschen2005distortion}, model order reduction (MOR)\cite{burkardt2006centroidal}\cite{du1999centroidal}, optimization \cite{okabe2000spatial}, and many other disparate fields.  The objective of this work is to present a family of multi-use, data-driven partitioning algorithms that concisely unifies many well-known methods with a single objective function and optimization routine.  The design is easy to use and interpret, and relies on few hyperparameters.  In addition, by alternating between local and global optimization objectives, our adaptive formulation is able to automatically uncover hidden structures within data.  As seen in Figure \ref{fig:spheres}, this can result in a significant improvement in performance over existing methodologies.  We hope to encourage cross-collaboration by demonstrating how problems from various application areas can be solved by a single algorithmic entity.  The remainder of the paper is broadly divided into two sections.  In the first, we introduce our family of algorithms and discuss its features.  The section following presents numerical experiments, and is itself divided into three sections, each devoted to a specific application area (subspace clustering, model order reduction (MOR), and matrix approximation).  After these sections, we close with final remarks and ideas for future work.\\

\noindent \textit{Notation}  Before continuing, we introduce some helpful notation.  At times we will use Matlab \cite{MATLAB} notation; e.g., if $A \in \mathbb{R}^{m \times n}$ and $j$ is a positive integer, $A(:,j)$ and $A(j,:)$ denote the $j^{th}$ column and row of $A$ respectively.  
For a collection of matrices, $\{B_i\}_{i=1}^k$, with $B_i \in \mathbb{R}^{m \times n}$, we let $\mathsf{diag}(B_i)$ represent the $km \times kn$ matrix
$$\mathsf{diag}(B_i) = \left ( \begin{array}{ccc}
                      B_1& & \\
                       &\ddots&\\
                       & &B_k\\
                      \end{array} \right ).$$
Given a set $B \subset \mathbb{R}^m$ and a point $x \in \mathbb{R}^m$, we denote by $B - x$ the set 
$\{ b - x \; | \; b \in B\}.$  If instead $B\in \mathbb{R}^{m \times n}$ is a matrix and $x \in \mathbb{R}^m$, the matrix $B - x \in \mathbb{R}^{m \times n}$ is given by
$B - x = [ B(:,1) - x \dots B(:,n) - x ].$ For the identity matrix in $\mathbb{R}^{m \times m}$, we write $I_m$.  We write $|V|$ to denote the cardinality of a set $V$, and for a natural number, $n$, $[n]$  denotes the set $\{1,\ldots,n\}$.

\section{A Parameter-Dependent Family of Algorithms}\label{section1}

We begin by outlining the algorithm family and its solution method.  This is followed by a description of the different objective function terms and indexing parameter.

Let $\Omega = \{x_i\}_{i=1}^n \subset \mathbb{R}^m$ be a data set consisting of $n$ points.  The general clustering objective is to partition $\Omega$ into $k$ disjoint sets (hereafter referred to as Voronoi sets), $V_i$, such that
$$\Omega = \bigcup_{i=1}^k V_i, \quad V_i \cap V_j = \emptyset \mbox{ for } i \neq j.$$
\noindent We assume $k$ is given for now.  Our algorithmic family is represented by the following functional equation, $\mathcal{G}_{\alpha}\left ( \{V_i\}_{i=1}^k , \{\Phi\}_{i=1}^k, \{\beta_i\}_{i=1}^k \right )$, indexed by the parameter $\alpha \in [0,1]$:

\begin{align*}
     \mathcal{G}_\alpha=&\sum_{i=1}^k\sum_{x \in V_i}\left ( \alpha\| x - \beta_i\|_2^2 + (1 - \alpha)\|(I - \Phi_i)(x - \beta_i)\|_2^2 \right )\\
    =&\sum_{i=1}^k\sum_{x \in V_i}\left ( \| x - \beta_i\|_2^2 - (1 - \alpha)\|\Phi_i(x - \beta_i)\|_2^2 \right ).\\
\end{align*}

\noindent Here, the $V_i$ partition $\Omega$, $\beta_i \in \mathbb{R}^m$, and $\Phi_i\in \mathbb{R}^{m \times m}$ are $d_i-$dimensional orthogonal projectors (hereafter referrred to as centroids).  Given $\Omega$, $k$, and a multi-index of projector dimensions, $\{d_i\}_{i=1}^k$, the optimization problem to be solved is:
$$ \min_{V_i,\Phi_i,\beta_i}\mathcal{G}_{\alpha}\quad \mbox{ such that }$$
$$\bigcup_{i=1}^k V_i = \Omega,\;\Phi_i^2 = \Phi_i,\; \mbox{rank}(\Phi_i) = d_i,\; i = 1,\ldots,k.$$
The solution is found via alternating minimization.  Given an initial partition with the $\beta_i$ set to $\frac{1}{|V_i|}\sum_{x \in V_i}x$, one repeats the following three steps until some convergence criteria is met:

\begin{enumerate}
    \item (Centroid Update) Hold the $V_i, \beta_i$ fixed and optimize over $\Phi_i$.
    \item (Voronoi Update) Fix $\beta_i, \Phi_i$ and optimize over $V_i$.
    \item (Mean Update)  Optimize over $\beta_i$ while keeping the $V_i$ and $\Phi_i$ fixed
\end{enumerate}

\noindent In the Centroid Update, we determine the optimal projectors.   These are given by $\Phi_i = U_iU_i^T,$ where $U_i \in \mathbb{R}^{m \times d_i}$ is the matrix whose columns contain the leading $d_i$ left singular vectors of $V_i - \beta_i$.  The Voronoi Update step proceeds by assigning points via the following rule:
$$ x \in V_i \mbox{ if }$$
$$i = \mbox{arg}\min_{l\in [k]}\left ( \|x - \beta_l\|_2^2 - (1-\alpha)\|U_l^T(x - \beta_l)\|_2^2\right )$$
\noindent Ties are broken by assigning points to the set with the smallest index.
For $\alpha \in \{0,1\}$, the optimal $\beta_i$ in the Mean Update step can be found via differentiation to be the Voronoi set means, $m_i = \frac{1}{|V_i|} \sum_{x \in V_i}x_i$.  For $\alpha \in (0,1)$, a closed-form solution is not available.  In this case, we perform a single gradient descent step to estimate $m_i$.  It can be shown that this process requires no line search and does not increase the overall complexity (See Section \ref{SIsection}). With these objects determined, we may rewrite our functional as:

\begin{equation}\label{functional}
   \mathcal{G}_{\alpha} =\sum_{i=1}^k\sum_{x \in V_i}\left ( \| x - m_i\|_2^2 - (1 - \alpha)\|U_i^T(x - m_i)\|_2^2 \right ). 
\end{equation}

\noindent In the simulations that follow, we halt this process when either (a) the number of iterations reaches 50 or (b), the difference between consecutive values of \ref{functional} falls below a user-supplied tolerance, $tol$ (e.g., $tol = 0.1$).  The value in (a) is selected based on many observations where the algorithm typically required fewer than ten iterations to satisfy $tol = 0.1$.  Since the alternating minimization process produces a non-increasing sequence of functional values, $\mathcal{G}_\alpha$ (see \cite{kerschen2002non} and \cite{du2003centroidal}), the tolerance criterion in (b) will ultimately be satisfied in the absence of an iteration limit.\\


\noindent \textit{Adaptation}  The optimization method just described keeps the number of Voronoi sets, $k$, and the projector dimensions, $\{d_i\}$, fixed.  In addition, the projectors are updated to solve a local minimization problem involving individual Voronoi sets:
$$\max_{\Phi_i\in \mathbb{R}^{m \times m}} \|\Phi_i(V_i - m_i)\|_F^2,\quad \mbox{ such that } $$
$$\Phi_i^2 = \Phi_i,\; \mbox{rank}(\Phi_i) = d_i,\; i \in [k].$$

\noindent We can convert this into a global optimization problem that involves all Voronoi sets by solving:
$$\max_{\Theta \in \mathbb{R}^{km \times km}}\|\Theta \mathsf{diag}(V_i - m_i)\|_F^2,\quad \mbox{such that }$$
$$\Theta^2 = \Theta,\; \mbox{rank}(\Theta) = r,\; r = \sum_{i=1}^k d_i.$$ 

\noindent The solution is given by $\mathsf{diag}(U_iU_i^T)$, where $U_i \in \mathbb{R}^{m \times \widetilde{d}_i}$ contains the leading $\widetilde{d}_i$ left singular vectors of $V_i - m_i$ that contribute to the dominant $r = \sum_{i=1}^kd_i$ dimensional subspace of $\mathsf{diag}(V_i - m_i)$.  It is this last characteristic that allows for adaptation.  For example, it may happen that no left singular vector from one or more of the $V_i-m_i$ contributes to the dominant $r$-dimensional subspace of $\mathsf{diag}(V_i-m_i)$.  In this case, the rank $r$ left singular matrix of $\mathsf{diag}(V_i-m_i)$ could look like the following:
$$\left ( \begin{array}{cccc}
          U_1U_1^T & & &\\
           & \ddots & &\\
           & & U_{k-1}U_{k-1}^T & \\
           & & & 0\\
           \end{array}\right )$$

\noindent When this occurs, we allow the number of sets, $k$, to change in order to match the number of singular matrices from each of the $V_i-m_i$ that contribute to the rank $r$ singular value decomposition (SVD) of $\mathsf{diag}(V_i-m_i)$. For the Voronoi sets that remain, we keep the corresponding $m_i$.  The result is a data driven routine in which the dimension and number of the sets $V_i$ are allowed to vary over the course of the algorithm.  Moreover, one can use the single value, $r$, as an input instead of specifying individual projector dimensions.\\  

\noindent We remark that the value of $\mathcal{G}_\alpha$ may increase with a reduction in $k$ (compare to the k-means objective value increasing with decreasing $k$ \cite{thorndike1953belongs}), albeit for one iteration.  However, as we now show, the method will still satisfy the tolerance criterion (b) mentioned above.  Given an initial value for $k$, the set of possible $k$ values encountered during the minimization process is finite and bounded below by one.  Moreover, the value of $k$ cannot increase.  Thus, at some point $k$ will become fixed, and the ensuing sequence of functional values, $\mathcal{G}_\alpha$, will be non-increasing.  Since $\mathcal{G}_\alpha \ge 0$, (b) must be satisfied. We can view this as the algorithm seeking to avoid local minima associated with an incorrect value of $k$, something that the functional \ref{functional} cannot capture. This also explains our three-step minimization process and its ordering.  For example, it is possible that some points $x_i$ will no longer have an assigned Voronoi set after adaptation.  This is remedied by retaining the means associated with the adapted centroids and performing the Voronoi Step prior to updating the $m_i$. As a result, we recommend supplying overestimates of $k$ and $r$ at initialization.\\

\noindent \textit{Choice of Indexing Parameter, $\alpha$.}
\noindent Here, we discuss the role of the indexing parameter, $\alpha$, and its implications; numerical examples appear in the following section.  As mentioned earlier, several well-known partitioning algorithms coincide with specific parameter settings of our family.  For example, setting $\alpha =1$ and $k>1$ results in the k-means algorithm \cite{macqueen1967some}.  With $\alpha = 0$, one is able to recover the k-subspaces (KSS)\footnote{As we show later, this objective also coincides with the Centroidal Voronoi Orthogonal Decomposition (CVOD) method \cite{du2003centroidal}} \cite{agarwal2004k} algorithm for subspace clustering by fixing the Voronoi set means to zero.  If instead the means are allowed to vary, one arrives at the Vector Quantization Principal Component Analysis (VQPCA) \cite{kerschen2002non}.  Reducing $k$ to one in this case results in the standard Principal Component Analysis (PCA) objective \cite{jolliffe2005principal}.  For $\alpha \in (0,1)$, the functional is a mixture of k-means and KSS-type terms.  One can view this case as a k-means objective with a regularization term based on Voronoi set subspace information.  Note that, for all parameter values, one can choose between an adaptive and non-adaptive setting (e.g., a standard or adaptive k-means/KSS algorithm). \\ 

\noindent \textit{Complexity and Hyperparameters.}  The algorithmic family as a whole roughly scales as $\mathcal{O}(rmn)$ in terms of complexity, and has few hyperparameters.  Excluding the choice of $\alpha$, the non-adaptive variants use the same number of parameters as KSS ($k$, a multi-index, $\{d_i\}_{i=1}^k$, and stopping tolerance, $tol$) while, as mentioned earlier, the adaptive version uses less ($k$,$r$,$tol$). Moreover, we have observed instances where the adaptive variants are able to automatically infer the correct parameter values (See Figure \ref{fig:spheres} and Section \ref{numerics}).

\section{Numerical Experiments}\label{numerics}

The purpose of this section is to demonstrate the performance of our algorithmic family with particular emphasis on the adaptive variants.  In addition, we hope to encourage the sharing of theory and algorithms between different areas by showing via our family that they often rely on similar tools.  The experiments are divided into three sections, with each focused on a specific application area.  These include subspace clustering, model order reduction, and matrix approximation.  Included within each section is a brief overview of the area and problem instance.

\subsection{Subspace Clustering}

In this section, we use our algorithms on a subspace clustering (SC) problem.  The data in an SC context is assumed to be from the union of subspaces, and the task is to (a) determine the number of subspaces (which may or may not overlap), (b) their dimensions, and (c) partition to the data points according to the subspace to which they are closest \cite{vidal2011subspace} \cite{wang2022convergence}.  Applications range from computer vision, gene analysis, image segmentation, temporal video segmentation, and image compression \cite{elhamifar2013sparse}.  Given the difficulty of the problem and its many applications, it is not surprising that many SC algorithms exist.  Our focus will be on the k-subspaces (KSS) algorithm \cite{agarwal2004k} \cite{tseng2000nearest}.  For more on the different (linear) SC algorithms and how they compare, the reader is referred to \cite{vidal2011subspace}; for extensions to nonlinear subspaces, see \cite{abdolali2021beyond}.

The KSS algorithm is a generalization of the well-known k-means \cite{hart2000pattern} routine, and corresponds to setting $\alpha = 0$ in our functional \ref{functional}:

$$\min_{\{(V_i,\Phi_i,\beta_i)\}_{i=1}^k} \sum_{i=1}^k \sum_{x \in V_i}||\beta_i - (I_m - \Phi_i)x||_2^2\quad \mbox{such that}$$
$$\Phi_i^2 = \Phi_i,\; \mbox{rank}(\Phi_i) = d_i,\; \sum_{i=1}^k d_i = r,$$ $$\beta_i \in \mathbb{R}^m,\quad \bigcup_{i=1}^k V_i = \Omega,\quad i = 1,\ldots,k.$$

\noindent In the case of non-affine subspaces, one can set the $\beta_i = 0$ \cite{agarwal2004k}.  The method is popular because of its simplicity and ease of use.  Moreover, its linear complexity in the number of points per iteration means it scales well to large problem instances \cite{wang2022convergence}.

\begin{figure}[H]
    \centering
    \begin{subfigure}[t]{\textwidth}
        \centering
        \includegraphics[width=1\textwidth,height=0.3\textwidth]{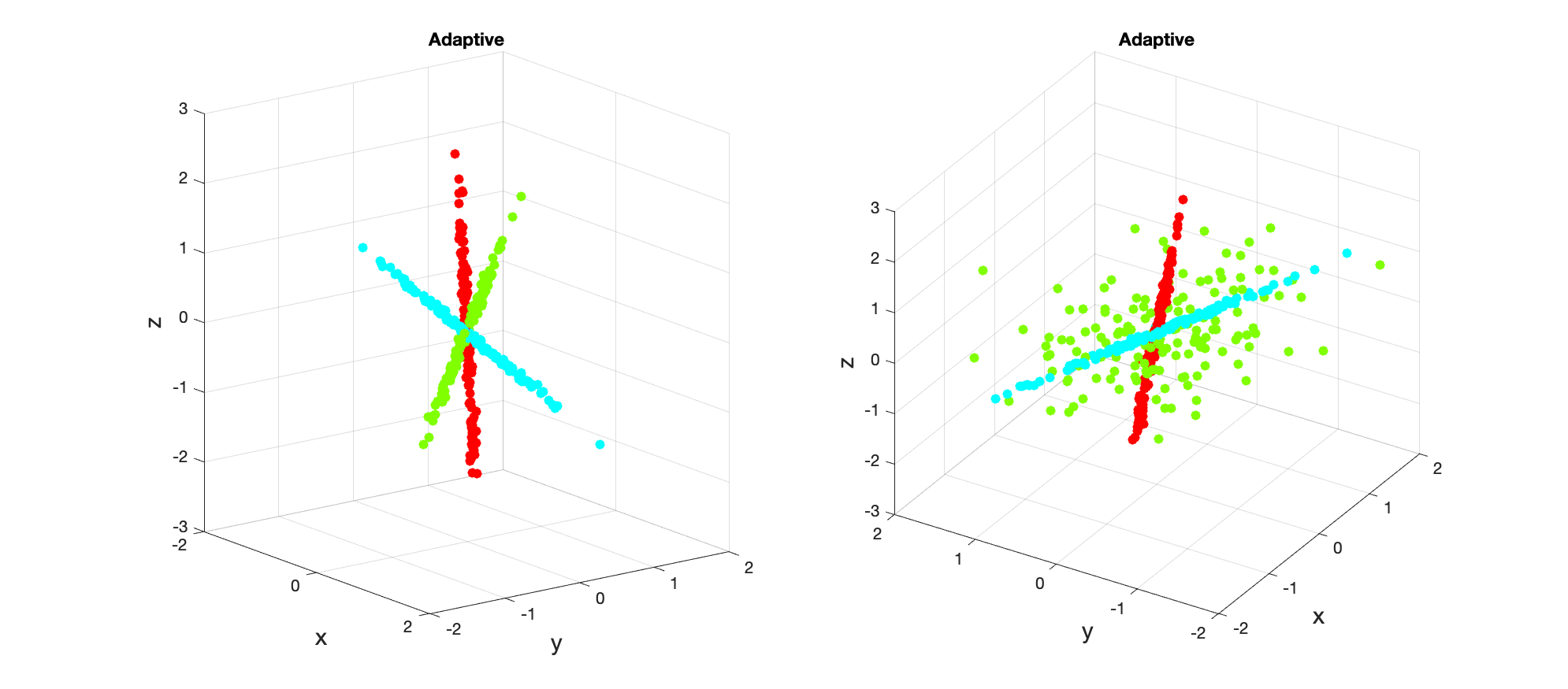}
        \caption{Adaptive KSS clustering result and rotated view.}
        \label{fig:sc1}
    \end{subfigure}
    \begin{subfigure}[t]{\textwidth}
        \centering
        \includegraphics[width=1\textwidth,height=0.3\textwidth]{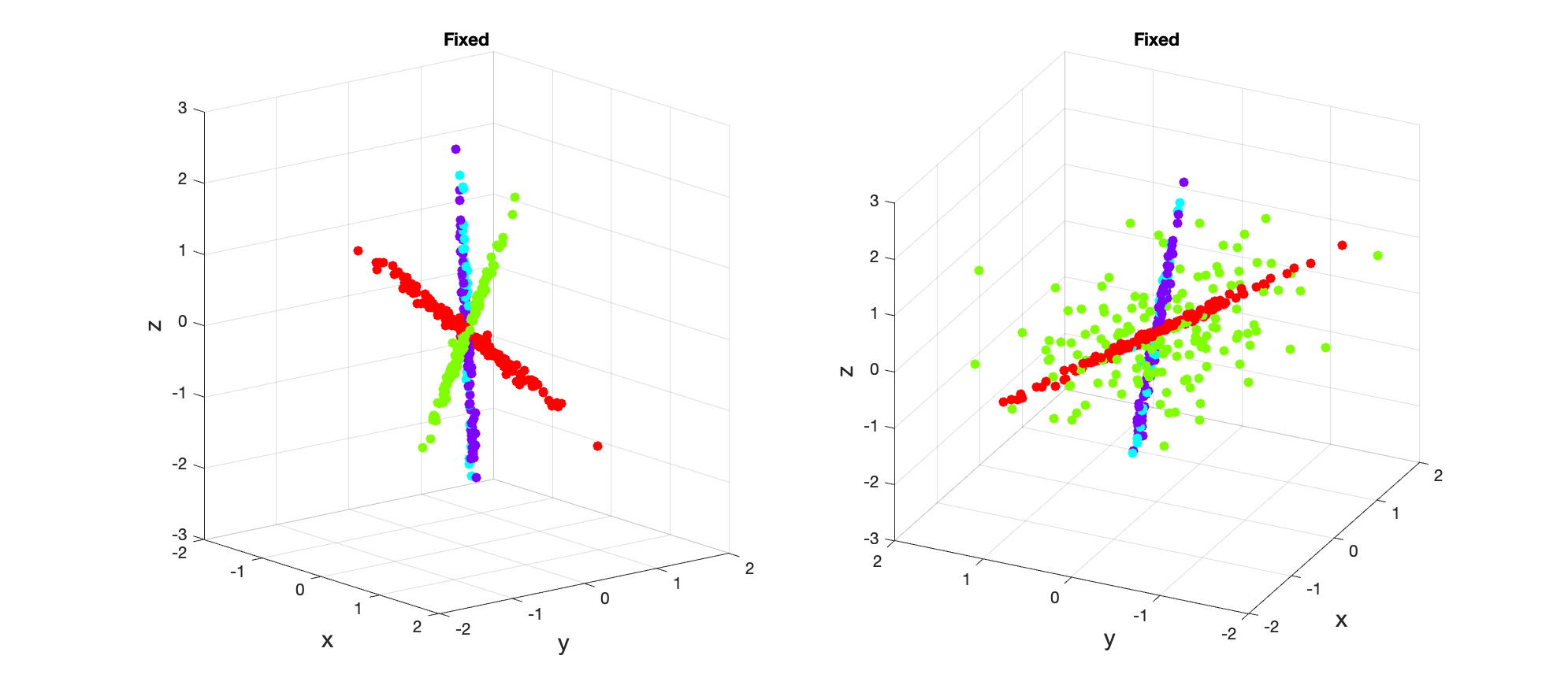}
        \caption{Non-adaptive KSS clustering result and rotated view using the same initial parameters as \ref{fig:sc1}.  With incorrect subspace dimensions and $k$, the algorithm is unable to capture the true clustering.}
         \label{fig:sc2}
    \end{subfigure}
    \caption{Subspace clustering result (colors) using the adaptive family with $\alpha = 0$, $k = 4$, and total dimension, $d_{total} = 7$.  Note that the algorithm is able to uncover the correct clusters and dimensions despite being initialized differently.}
    \label{fig:sc}
\end{figure}

Figures \ref{fig:sc1} and \ref{fig:sc2} show a simple example that showcases our algorithm's adaptive mechanism in an SC setting.  The data consist of two intersecting planes along with an intersecting line in $\mathbb{R}^3$.  Thus, the true subspace dimensions are $(2,2,1)$ with $k_{true} = 3$ subspaces.  We use our adaptive algorithm variant with $\alpha = 0$ with no restrictions on the $m_i$.  In addition, we initialize the total subspace dimension to seven and the number of subspaces to four.  To address the KSS method's known sensitivity to initial partitions, we employ the ensemble approach developed in \cite{lipor2021subspace}.  The basic idea is to run the algorithm $B$ times (e.g., $B = 1000$, though we use 200).  One defines the matrix $A \in \mathbb{R}^{n \times n}$ where $a_{ij}$ gives the average number of instances where the points $x_i$ and $x_j$ are co-clustered.  After thresholding each row and column of $A$ so that each row and column only contains the top $q$ entries, one performs spectral clustering \cite{ng2001spectral} \cite{von2007tutorial}.  For our case, we set $q = 40$ arbitrarily, though there are more-informed choices \cite{heckel2014neighborhood} \cite{heckel2015robust}.  We also observe good results without thresholding, though this is likely a problem-dependent artifact.

The spectral clustering step discussed above requires one to specify a value for the number of clusters (i.e., for k-means).  Although one may repeat the k-means step for several values of $k$, and in some instances one has an idea for the correct $k$ \cite{lipor2021subspace}, we choose to allow our algorithm to provide an estimate.  In particular, we run two versions of our algorithm for each of the $B = 200$ iterations on the same data set above.  We use $\alpha = 0.5$ and a stopping tolerance, $ tol = 0.1$, for the first, and $(\alpha = 0,\;tol = 0.01)$ for the second.  We choose $\alpha = 0.5$ for demonstration purposes in an attempt to balance the KSS and k-means terms in \ref{functional}. Optimizing the choice of this parameter will be left to future work. Both algorithms are allowed to adapt. The idea is to use the first instance to provide an estimate for $k$ (and not a perfect partitioning, hence the larger value of $tol$), and the second to focus on subspace estimation.  Once complete, we use the average final $k$ value from the $\alpha = 0.5$ runs along with the affinity matrix from the $\alpha = 0$ runs as input to a spectral clustering routine. 

The results are shown in Figure \ref{fig:sc}, along with that from an ensemble KSS (non-adaptive) implementation with inputs $k = 4$, $d = (2,2,2,1)$ (i.e., the initial $k$ and total dimensions match that of the adaptive runs).  For the KSS implementation, we use $k = 4$ in the spectral clustering step as well, with the intent of (a) demonstrating the effect of not knowing the correct $k$ and (b) highlighting our adaptive mechanism.

Figure \ref{fig:sc} suggests that our adaptation approach has the ability to determine the appropriate subspace parameters and partitioning, even when provided incorrect starting values.  In fact, the $\alpha = 0.5$ instance determined the correct $k$ $86\%$ of the time.  Note similar behavior is also observed in the case of high-dimensional spheres shown in Figure \ref{fig:spheres}, but for a generic clustering task.

We remark that the articles \cite{lipor2021subspace} and \cite{lane2019adaptive} also present adaptive KSS-type algorithms as well as unique initialization strategies.  In the first work, the dimensions for each Voronoi set are updated based on the index of the largest eigengap from the covariance matrix constructed from $x \in V_i$.  However, the number of sets $V_i$ remains fixed.  Updates in the second work are achieved via several regularization terms in addition to the usual KSS objective.  We note, however, that the subspace basis elements are not required to be orthogonal.  In contrast to our work and \cite{lipor2021subspace}, the method in \cite{lane2019adaptive} requires a larger number of hyper-parameters (e.g., for initialization, stochastic gradient descent, and the algorithm itself).

\subsection{Model Order Reduction}\label{MORsection}

This section focuses on applying our algorithmic family to problems related to Model Order Reduction (MOR).  In a typical scenario from MOR, one is given a large system of high-dimensional equations; e.g., a system of ordinary differential equations (ODEs) resulting from a discretization of a system of partial differential equations (PDE) \cite{alla2017nonlinear}.  In an effort to reduce computational complexity and runtime, one can project this system onto a lower-dimensional space \cite{alla2022first} \cite{quarteroni2015reduced}.  In practice, one determines the low-dimensional basis via output from several (expensive) model runs with different parameter values.  Referred to as snapshots, these model output vectors can be collected into a single data matrix.  A popular approach called Proper Orthogonal Decomposition (POD), results from computing the SVD of the snapshot matrix, and forming a basis from a subset of the leading left singular vectors.  A lower-dimensional system is achieved by projecting the system equations onto this reduced basis in a Galerkin-type fashion \cite{quarteroni2010numerical}.  For more on this topic, see \cite{quarteroni2010numerical} \cite{quarteroni2015reduced} and the references therein.

While POD remains a popular MOR method, other approaches exist.  These include Dynamic Mode Decomposition (DMD) \cite{alla2017nonlinear}\cite{schmid2022dynamic} \cite{tu2013dynamic} and Discrete Empirical Interpolation Method (DEIM) \cite{chaturantabut2010nonlinear}, as well as methods that partition the snapshot set prior to forming a low-dimensional basis.  This last refers to Centroidal Voronoi Tessellation (CVT) type methods \cite{burkardt2006centroidal}\cite{okabe2000spatial}, which include Vector Quantization Principal Component Analysis (VQPCA)\cite{kerschen2002non} \cite{kerschen2005distortion} and the Centroidal Voronoi Orthogonal Decomposition (CVOD) \cite{du1999centroidal} method.  

Given a data set $\{x_i\}_{i=1}^n = \Omega \subset \mathbb{R}^m$, VQPCA attempts to solve the following:

$$\min_{\{(V_i,\Phi_i,\beta_i)\}_{i=1}^k} \sum_{i=1}^k \sum_{x \in V_i}||\beta_i - (I_m - \Phi_i)x||_2^2\quad \mbox{such that}$$
$$\Phi_i^2 = \Phi_i,\; \mbox{rank}(\Phi_i) = d_i,\; \sum_{i=1}^k d_i = r,$$ $$\beta_i \in \mathbb{R}^m,\quad \bigcup_{i=1}^k V_i = \Omega,\quad i = 1,\ldots,k.$$
\noindent CVOD minimizes a similar objective obtained by removing the $\beta_i$:

$$\min_{\{(V_i,\Phi_i)\}_{i=1}^k} \sum_{i=1}^k \sum_{x \in V_i}||(I_m - \Phi_i)x||_2^2\quad \mbox{such that}$$
$$\Phi_i^2 = \Phi_i,\; \mbox{rank}(\Phi_i) = d_i,\; \sum_{i=1}^k d_i = r,\quad \bigcup_{i=1}^k V_i = \Omega,\quad i \in [k].$$

\noindent In both cases, $d_i$ denotes the low-dimensional basis for each $V_i$ and $r = \sum_{i=1}^kd_i$ is the total dimension.  In terms of our algorithmic family, these algorithms both correspond to $\alpha = 0$ while the $m_i = 0$ specifically recovers the CVOD routine.  In addition, we see that these objectives are exactly the same as the KSS method discussed in the previous section.

We use an example problem from \cite{alla2022first} to demonstrate the performance of our adaptive algorithm in an MOR setting.  The system there is given by:
\begin{eqnarray*}
\frac{\partial z}{\partial t}(w,t)&=& \alpha \left( \frac{\partial^2}{\partial x^2} z(w,t) + \frac{\partial^2}{\partial y^2} z(w,t)\right ) + f(t,z(w,t))\\
z(w,0)&=& z_0(w)\\
z(w,t) &=& 0,\; (w,t) \in \partial \Omega \times [0,T]\\
\end{eqnarray*}

\noindent We convert this into a system of ODEs by using a finite difference approximation with spatial steps $\Delta x = \Delta y = 0.0125$:

\begin{eqnarray}
M\dot{z}(w,t) &=& Bz(w,t) + f(t,z(w,t)),\quad (w,t) \in \Omega \times [0,T] \nonumber \\ 
z(0) &=& z_0(w) \nonumber \\
z(w,t)&=& 0,\quad (w,t) \in \partial \Omega \times [0,T]\nonumber
\end{eqnarray}
where $z_0 = (x_0,y_0) \in \Omega \subset \mathbb{R}^2$ is the initial condition, $t \in [0, T]$, and $M,B \in \mathbb{R}^{n \times m}$ are given matrices.  These types of systems appear when discretizing PDEs for heat transfer or wave equations \cite{alla2022first} \cite{quarteroni2010numerical}.  We take $M$ as the identity, $B$ is a discrete Laplacian operator formed via finite differences, $f(t,z(t)) = 10\left ( z(t)^2 - z(t)^3 \right )$, $\alpha = 0.05$, and $z_0(w) = \sin(\pi x) \sin(\pi y)$. Our solution uses an implicit Euler scheme with Newton's method and the following parameters: $\Delta T = 0.05$, $T = 2$, and $\Omega = [0 ,1].$  The dimension of the solution space is $\mathbb{R}^{6241}$.


\begin{figure}[H]
    \centering
    \begin{subfigure}[t]{0.45\textwidth}
        \centering
        \includegraphics[width=1\textwidth,height=0.85\textwidth]{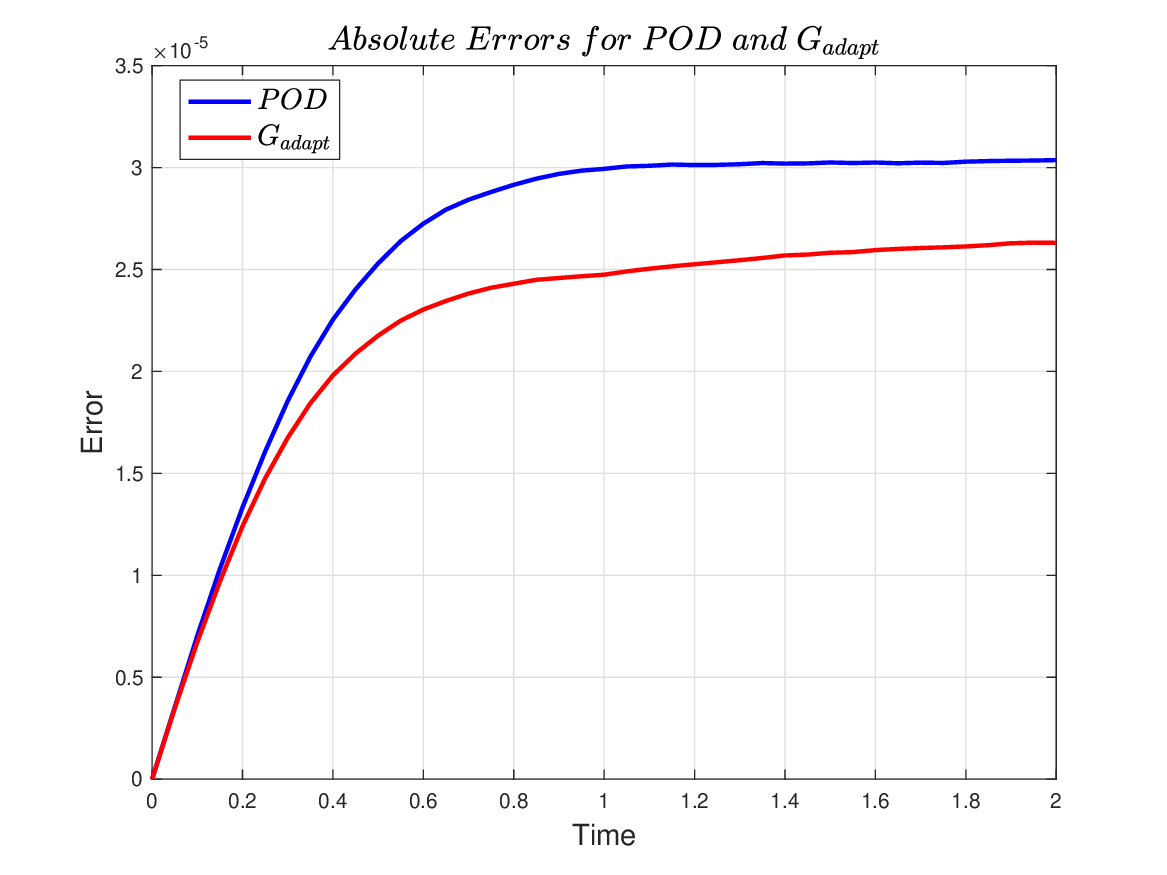}
        \caption{}
        \label{fig:morfigure1}
    \end{subfigure}
    \hfill
    \begin{subfigure}[t]{0.5\textwidth}
        \centering
        \includegraphics[width=1\textwidth,height=0.85\textwidth]{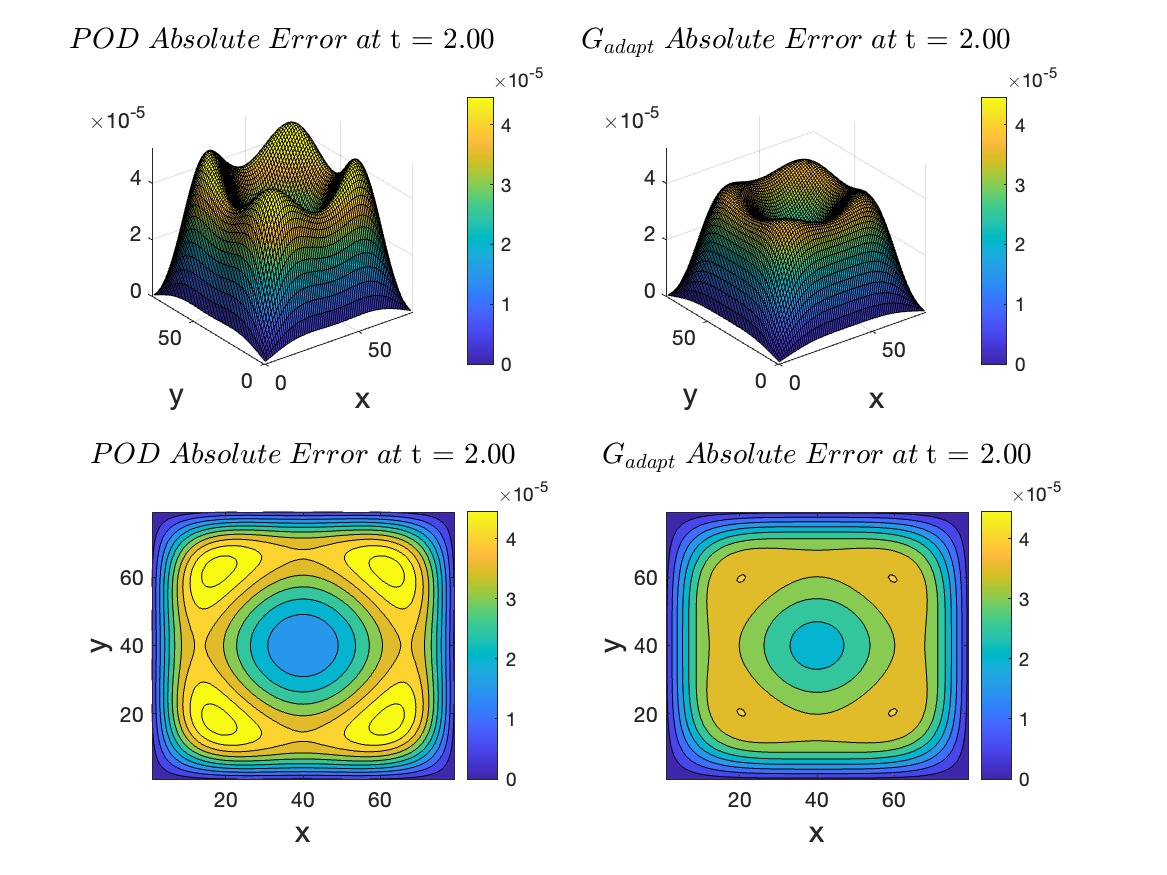}
        \caption{}
         \label{fig:morfigure2}
    \end{subfigure}
    \caption{Error comparison for POD and our adaptive family (labeled $G_{adapt}$) in time (Figure \ref{fig:morfigure1}) and at time $t = 2$ (Figure \ref{fig:morfigure2}).  See text for algorithm settings.}
    \label{fig:morfigure}
\end{figure}

Figure \ref{fig:morfigure} shows the results comparing POD to our adaptive variant with $\alpha =0$ and $m_i = 0$.  Both methods are set to return ten basis elements.  We use $k_{init} = 10$ and $tol = 10^{-4}$.  Instead of initializing the Voronoi sets randomly, we use an approach based on ranking the correlations between each pair of snapshots (See Section \ref{SIsection} for details).  The figure suggests that partitioning combined with an informed choice of initialization (an idea echoed in the previous section) can lead to a significant improvement in performance.  For more on this algorithm in the context of MOR, see \cite{emelianenko2024optimality}.

\subsection{Matrix Approximation}

This section introduces two related matrix approximation problems and shows how our family can be used to approximate their solution.  The first is referred to as the column-subset selection problem (CSSP) \cite{boutsidis2009improved}.  Given a matrix $A \in \mathbb{R}^{m \times n}$ with $\mbox{rank}(A) = \rho$, and a target rank, $0<r \le \rho,$ the goal of CSSP is to form $C \in \mathbb{R}^{m \times r}$ consisting of $r$ columns of $A$ that minimizes $$||(I - CC^\dagger)A||_\xi,\quad \xi \in \{2,F\},$$ over all possible $m \times r$ matrices $C$ whose columns are taken from $A$ ($C^\dagger$ denotes the Moore-Penrose pseudoinverse of the matrix $C$).  This problem is difficult to solve \cite{shitov2017column}, since determining the best solution requires enumerating all $\binom{n}{r}$ possible solutions.  Solution methods range from deterministic approaches, where one applies a classical matrix factorization routine to select columns (e.g., LU-decomposition, QR-decomposition with partial pivoting), to probabilistic techniques \cite{dong2021simpler}.  This last group of methods, where columns are selected via carefully made probability distributions, have been shown to scale well with large problem sizes and lead to improved upper and lower error bounds \cite{mahoney2009cur} \cite{derezinski2021determinantal}.  

\begin{figure}[H]
    \centering
    \begin{subfigure}[t]{0.45\textwidth}
        \centering
        \includegraphics[width=\textwidth,height=0.8\textwidth]{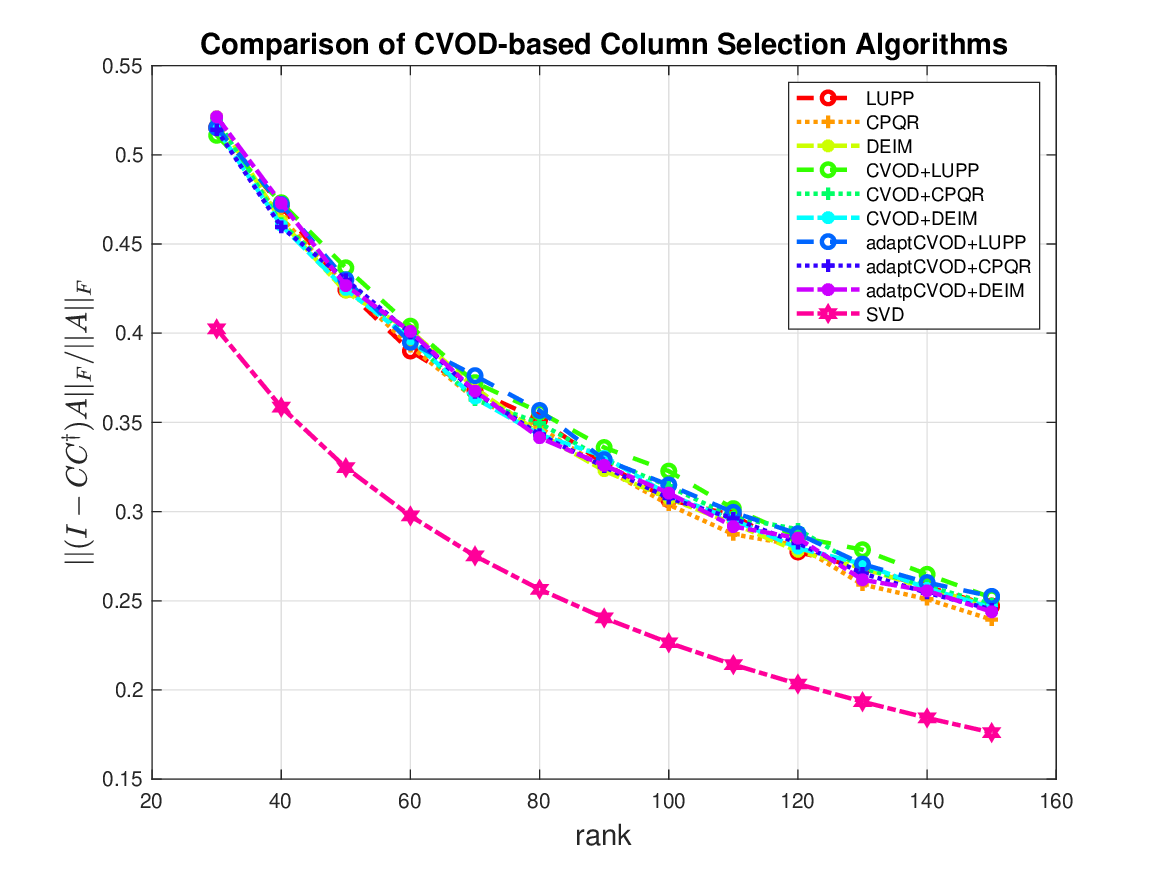}
        \caption{Results for our $(\alpha = 0,m_i = 0)$ variants.}
        \label{fig:matrix1}
    \end{subfigure}
    \begin{subfigure}[t]{0.45\textwidth}
        \centering
        \includegraphics[width=\textwidth,height=0.8\textwidth]{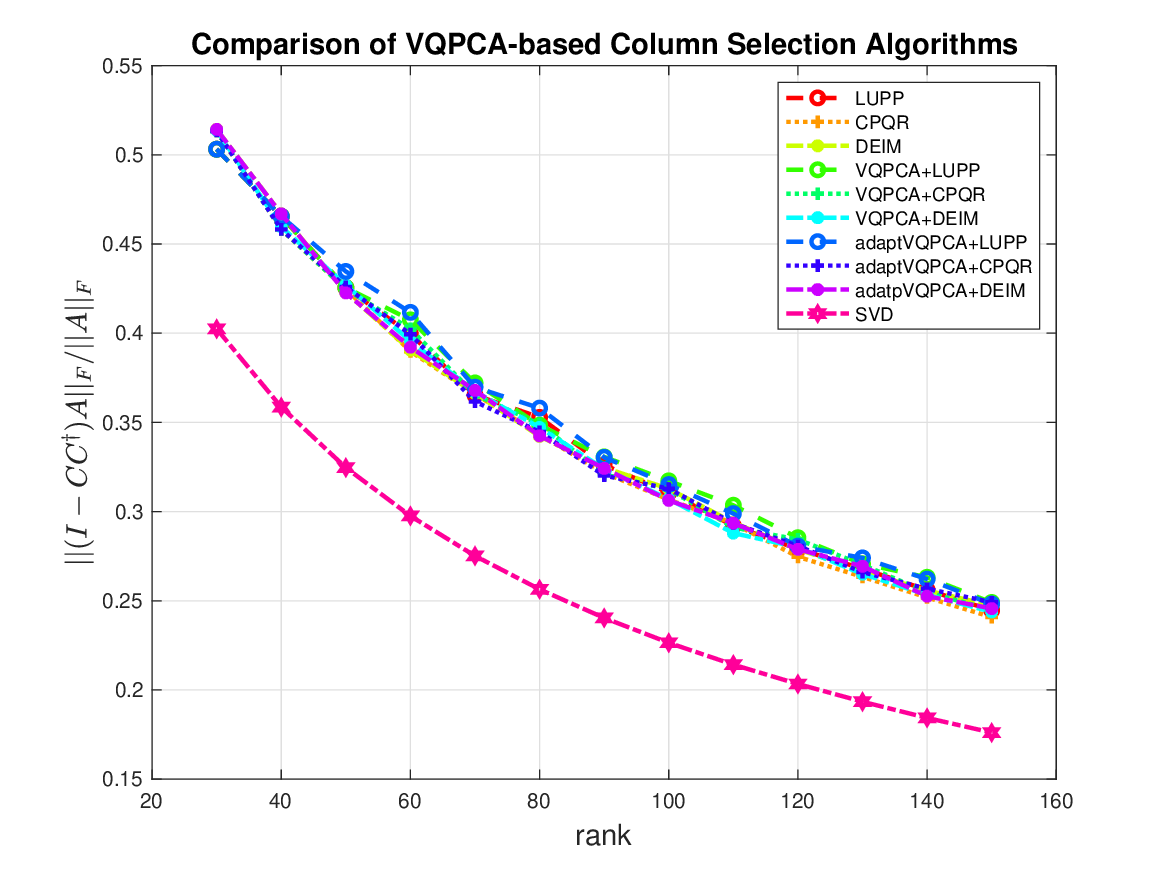}
        \caption{Same as \ref{fig:matrix1} but for $\alpha = 0$ and variable $m_i$.}
        \label{fig:matrix2}
    \end{subfigure}
    \caption{Matrix approximation error results for various algorithms from our family.  The data consists of the matrix $A \in \mathbb{R}^{60000 \times 785}$ containing MNIST training images.  Images taken from \cite{emelianenko2024optimality}.}
    \label{fig:matrix}
\end{figure}

These same column-selection techniques can be used to construct CUR decompositions \cite{goreinov1997theory} \cite{hamm2019cur} \cite{voronin2017efficient}\cite{wang2013improving}.  In this setting, a matrix $A \in \mathbb{R}^{m \times n}$ is factored into the product of three matrices, $A \approx CUR$, where $C\in \mathbb{R}^{m \times r}, R \in \mathbb{R}^{r \times n}$ contain, respectively, columns and rows selected from $A$.  The matrix $U \in \mathbb{R}^{r \times r}$ is chosen to make the residual, $A - CUR$, small.  Options include $U = C^\dagger A R^\dagger$ and $U = A(I,J)^\dagger$, where $I,J$ denote the row and columns indices used to form $C$ and $R$.  One advantage of the latter version is that one does not have to revisit the original matrix \cite{hamm2020stability}.  These factorizations are attractive since, unlike an SVD, they preserve attributes from the parent matrix (e.g., sparsity, non-negativity, interpretability).  Applications include recommendation systems analysis, DNA analysis, and hyper-spectral image analysis \cite{mahoney2009cur}.

Our algorithmic family works in this setting by partitioning the columns/rows of a matrix and then applying an established CSSP algorithm to each piece.  As theorem \ref{CSSPtheorem} from \cite{emelianenko2024optimality} shows, this results in reconstruction errors whose upper bounds are dominated by the quality of the final partition.

\begin{theorem}\label{CSSPtheorem}
   Let $A \in \mathbb{R}^{m \times n}$, rank($A$) = $\rho$, and $0<r \le \rho$, $0<k< n$ be integers.  If $C \in \mathbb{R}^{m \times r}$ is the output from a CSSP routine paired with an adaptive or non-adaptive variant from the family defined by $\mathcal{G}_\alpha$ with $\alpha = 0$ and $m_i = 0$,  then
   \begin{align*}
       \|(I_m - CC^\dagger)A\|_F &\sim \mathcal{O}  \left (\mathcal{G}^* \right )
   \end{align*}
    where $\mathcal{G}^*$ is the energy value of either of these variants at completion.
\end{theorem}

For more details and theory, see \cite{emelianenko2024optimality}.  As an example, we take $A \in \mathbb{R}^{60000\times 784}$ where the rows consist of the MNIST training image data.  For ranks $r = 30,40,\ldots,150$ and $k_{init} = 5$, we record the reconstruction error, $\frac{\|(I - CC^\dagger)A\|_F^2}{\|A\|_F^2}$.  The algorithm variants we use all correspond to the $\alpha = 0$ case.  Versions where the means are fixed to zero are termed CVOD-based, and those with variable means are referred to as VQPCA-based.  We pair (and compare) our variants to the following well-known CSSP methods:  LU factorization with partial pivoting (LUPP), Column-pivoted QR decomposition (CPQR), and DEIM \cite{sorensen2016deim}.  For more about these, see \cite{dong2021simpler}.

The results in Figure \ref{fig:matrix} show our variants to be competitive with existing approaches.  Moreover, a feasible per-iteration complexity and potential for parallelization mean these methods can scale to large problem sizes.  We remark that much work has been done in regards to reducing the complexity of CSSP algorithms while improving or maintaining error guarantees \cite{derezinski2020exact}\cite{gidisu2021hybrid}.  With this in mind, we believe that our methods would be well-suited as a way to supplement existing approaches that can make use of clustering in addition to a CSSP solution.  For example, the partitioning ability of our variants may help reduce the complexity involved in approximating parameter-dependent matrices \cite{park2024low}.

\section{Conclusion}
In this paper, we present a parametrized suite of data-adaptive partitioning algorithms that unites and enhances several well-known methods. Our experiments using problems from different fields highlight the collection's versatility and demonstrate the ability of our adaptive process to automatically determine hidden data structures and patterns in large, high-dimensional datasets.  In addition, these simulations suggest the potential for pushing the boundaries of existing methodologies by fusing ideas from different fields.  For example, the ensemble method and theoretical results found in \cite{lipor2021subspace} and \cite{wang2022convergence} may improve the partitioned-based MOR methods discussed in section \ref{MORsection} as well as have applications to dynamical systems.  Our candidates for future work include investigating the role of $\alpha$ in our functional \ref{functional} and how its inclusion is sensitive to the choice of problem.  For example, setting $\alpha = 0$ in a general (not SC) clustering task typically yields poor results, while $\alpha \in (0,1)$ does well.  On the other hand, setting $\alpha \neq 0$ for SC tasks nearly always results in poor performance.  Another avenue is to see if runs using multiple $\alpha$ values can inform the choice of the thresholding parameter, $q$, in the SC ensemble method outlined in \cite{lipor2021subspace}.  A more ambitious task would be to align our adaptive mechanism with the theoretical results from \cite{lipor2021subspace} and \cite{wang2022convergence}. 

\bibliographystyle{unsrt}
\bibliography{main}

\newpage

\section{Appendix}\label{SIsection}
The purpose of this supplement is to provide details on (1) the gradient descent formulation used used within our minimization routine and (2), the initialization scheme used in our model order reduction (MOR) approach.

\subsection{Gradient Descent}

In this section we investigate the optimization procedure for our functional, $\mathcal{G}$, for $ \alpha \in (0,1)$; see below.

\begin{eqnarray*}
    \mathcal{G}_{\alpha}\left ( \{V_i\}_{i=1}^k , \{\Phi\}_{i=1}^k, \{m_i\}_{i=1}^k \right ) &=&\sum_{i=1}^k\sum_{x \in V_i}\left ( \alpha\| x - m_i\|_2^2 + (1 - \alpha)\|(I - \Phi_i)(x - m_i)\|_2^2 \right )\\
    &=&\sum_{i=1}^k\sum_{x \in V_i}\left ( \| x - m_i\|_2^2 - (1 - \alpha)\|\Phi_i(x - m_i)\|_2^2 \right ).\\
\end{eqnarray*}

For $\alpha \in \{0,1\}$, the summand can be written as a single norm, thus making optimization for the $m_i$ straightforward (i.e., one can find closed form solutions).  The case for $\alpha \in (0,1)$ is not as clear.  In this section, we derive a gradient descent solution for determining the $m_i$ that does not rely on a line search.  

We begin by determining the derivatives of the following with respect to $m \in \mathbb{R}^m$:
$$\|x - m\|_2^2 - (1 - \alpha)\|\Phi(x - m)\|_2^2.$$
\noindent Recall that $x\in \mathbb{R}^{m}$ and that $\Phi \in \mathbb{R}^{m \times m}$ is an orthogonal projector.
Starting with the first term on the left, we may expand:
$$\|x - m\|_2^2 = \|x\|_2^2 - 2\langle x, m \rangle + \|m\|_2^2.$$

\noindent By taking the derivative with respect to an arbitrary component, $m_i$, one can show
$$\frac{\partial }{\partial m}\langle x,m \rangle = x,\quad \frac{\partial }{\partial m}\|m\|_2^2 = 2m.$$
Combining these gives
$$\frac{\partial }{\partial m}\|x - m\|_2^2 = 2(m - x).$$

\noindent Taking the left-most term, we may expand:
$$\|\Phi(x - m)\|_2^2 = \|\Phi x\|_2^2 - 2\langle \Phi x, \Phi m \rangle + \|\Phi m\|_2^2.$$

\noindent Keeping in mind that $\Phi$ is an orthogonal projector, one can show 
$$\frac{\partial }{\partial m} \langle \Phi x, \Phi m \rangle = \Phi x,\quad \frac{\partial }{\partial m} \|\Phi m\|_2^2 = 2 \Phi m.$$

\noindent Combining our results gives
$$\frac{\partial }{\partial m}( \|x - m\|_2^2 - (1 - \alpha)\|\Phi(x - m)\|_2^2) = 2(m - x) - 2(1 - \alpha)\Phi m.$$

We now extend this result to the full functional, $\mathcal{G}$, and each $m_i \in \mathbb{R}^m.$

\noindent We have
\begin{eqnarray*}
    \frac{\partial }{\partial m_i}\mathcal{G} &=&\frac{\partial }{\partial m_i} \sum_{i=1}^k\sum_{x \in V_i}\left ( \| x - m_i\|_2^2 - (1 - \alpha)\|\Phi_i(x - m_i)\|_2^2 \right ) \\
    &=& \sum_{x \in V_i}( 2(m_i - x) - 2(1-\alpha)\Phi_i (m_i - x))\\
    &=&2 \sum_{x \in V_i}\left ( (m_i - x) - (1 - \alpha)\Phi_i(m_i-x) \right )\\
    &=& 2 \left ( (|V_i|m_i - |V_i|\bar{x}_i ) - (1 - \alpha)\Phi_i (|V_i|m_i - |V_i|\bar{x}_i ) \right )\\
    &=& 2 |V_i| \left ( (m_i - \bar{x}_i) - ( 1- \alpha)\Phi_i (m_i - \bar{x}_i) \right )\\
    &=& 2|V_i| (I - (1 - \alpha )\Phi_i)(m_i - \bar{x}_i)\\
    &=& 2|V_i|\Gamma_i (m_i - \bar{x}_i)\\
    &\equiv& y_i.\\
\end{eqnarray*}
\noindent Here, $\bar{x}_i = \frac{1}{|V_i|}\sum_{x \in V_i}x$ and $\Gamma _i = I - (1 - \alpha )\Phi_i.$
Note that $y_i$ represents the gradient direction of $\mathcal{G}$ with respect to $m_i$.

\subsubsection{Step Size}

Now we use the above analysis to determine the appropriate gradient step.  Recall that the typical gradient descent format is given by
$$ m_i^{(l+1)} = m_i^{(l)} - \gamma_i \frac{\partial }{\partial m_i}\mathcal{G}_\alpha = m_i^{(l)} - \gamma_i y_i,$$
\noindent where $(l)$ denotes the iteration counter.

\noindent Since our goal is to update the means, $m_i$, such that the local Voronoi energy decreases, we want the new $m_i^{(l+1)}$ to (ideally) satisfy:

\begin{equation}\label{inequality}
\sum_{x \in V_i}( \|x - m_i^{(l+1)}\|_2^2 - (1 - \alpha)\|\Phi_i(x - m_i^{(l+1)}\|_2^2) < \sum_{x \in V_i}( \|x - m_i^{(l)}\|_2^2 - (1 - \alpha)\|\Phi_i(x - m_i^{(l)}\|_2^2).
\end{equation}
Our approach will be to substitute $m_i^{(l)} - \gamma_i y_i$ for $m_i^{(l+1)}$ and determine a sufficient $\gamma_i$.

\noindent Starting with the left term in \ref{inequality}, we have

\begin{eqnarray*}
    \sum_{x \in V_i}( \|x - m_i^{(l+1)}\|_2^2 - (1 - \alpha)\|\Phi_i(x - m_i^{(l+1)}\|_2^2)&=& \sum_{x \in V_i}( \|x - m_i^{(l)} + \gamma_i y_i\|_2^2 - (1 - \alpha)\|\Phi_i(x - m_i^{(l)}+ \gamma_i y_i\|_2^2)\\
\end{eqnarray*}

$$ = \sum_{x \in V_i}\left [ \left (\|x - m_i^{(l)}\|_2 + 2\langle x - m_i^{(l)}, \gamma_i y_i \rangle + \|\gamma_i y_i\|_2^2 \right ) - (1 -\alpha) \left ( \|\Phi_i (x - m_i^{(l)})\|_2 + 2 \langle \Phi(x - m_i^{(l)}),\gamma \Phi_i y_i \rangle + \|\gamma_i \Phi_i y_i\|_2^2  \right ) \right ]$$\\

\noindent Note that the first terms from both groups in parentheses will cancel the entire term on the right hand side of \ref{inequality}.  After canceling, the inequality takes the form:

\begin{align*}
0&>\sum_{x \in V_i} \left ( 2\langle x - m_i^{(l)}, \gamma_i y_i \rangle + \|\gamma_i y_i\|_2^2 - 2(1 -\alpha) \langle \Phi(x - m_i^{(l)}),\gamma \Phi_i y_i \rangle - (1 - \alpha)\|\gamma_i \Phi_i y_i\|_2^2 \right )\\
&=\sum_{i \in V_i} \left ( 2\langle x - m_i^{(l)}, \gamma_i y_i \rangle- 2(1 -\alpha) \langle \Phi(x - m_i^{(l)}),\gamma \Phi_i y_i \rangle + \|\gamma_i y_i\|_2^2-(1 - \alpha)\|\gamma_i \Phi_i y_i\|_2^2 \right )\\
&=\sum_{i \in V_i} \left ( 2\gamma_i\langle x - m_i^{(l)},  y_i \rangle- 2\gamma_i \langle (1 -\alpha)\Phi(x - m_i^{(l)}), y_i \rangle + \gamma_i^2\| y_i\|_2^2-\gamma_i^2(1 - \alpha)\|\Phi_i y_i\|_2^2 \right )\\
&=\sum_{i \in V_i} \left ( 2\gamma_i\langle x - m_i^{(l)} - (1 - \alpha)\Phi_i(x - m_i^{(l)}),  y_i \rangle + \gamma_i^2 \left (\| y_i\|_2^2-(1 - \alpha)\|\Phi_i y_i\|_2^2 \right )\right )\\
\end{align*}
\begin{align*}
   \Rightarrow -2 \gamma_i \sum_{x \in V_i} \langle \Gamma_i(x - m_i^{(l)}),y_i \rangle & >  \sum_{x \in V_i }\gamma_i^2 \left (\| y_i\|_2^2-(1 - \alpha)\|\Phi_i y_i\|_2^2 \right )\\
    \Rightarrow -2 \gamma_i \sum_{x \in V_i} \langle \Gamma_i(x - m_i^{(l)}),y_i \rangle &> \gamma_i^2|V_i|\left (\| y_i\|_2^2-(1 - \alpha)\|\Phi_i y_i\|_2^2 \right )\\
     \Rightarrow -2 \sum_{x \in V_i} \langle \Gamma_i(x - m_i^{(l)}),y_i \rangle &> \gamma_i|V_i|\left (\| y_i\|_2^2-(1 - \alpha)\|\Phi_i y_i\|_2^2 \right )\\
\end{align*}

\noindent This last follows since $\gamma_i >0$.  Let
$$\xi_i = -2 \sum_{x \in V_i} \langle \Gamma_i(x - m_i^{(l)}),y_i \rangle, \quad \zeta_i = |V_i|\left (\| y_i\|_2^2-(1 - \alpha)\|\Phi_i y_i\|_2^2 \right ).$$

\noindent Then the inequality will be satisfied by 
$$ \gamma_i = \frac{\xi_i}{\zeta_i + \eta}$$
\noindent for any $\eta >0$.

\noindent In practice, one should take $\gamma_i = 0$ if $\xi_i <0$; i.e., simply don't update the mean.  Note that when $\xi_i >0$, there is no restriction on the value of $\eta$.  

\subsubsection{Complexity}

In this section, we determine complexity for computing the $\gamma_i$, $i = 1,\ldots,k$.  By expanding the inner product, $\langle \Gamma_i(x - m_i^{(l)}),y_i \rangle$ and using the definition of $y_i$, one can show 
$$\left \langle \Gamma_i(x - m_i^{(l)}),y_i \right \rangle = c_i \left \langle x - m^{(l)}_i, m^{(l)}_i - \bar{x} \right \rangle + f_i \left \langle U_i^T(x - m_i^{(l)}),U^T_i(m_i^{(l)}-\bar{x}) \right \rangle,$$
\noindent where $c_i = -2|V_i|$ and $f_i = \left ( (1-\alpha)^2 - (1 - \alpha) - c_i(1 - \alpha) \right ).$  With $m_i^{(l)},x, \bar{x} \in \mathbb{R}^m$ and $U_i \in \mathbb{R}^{m \times d_i}$, the complexity of computing this inner product is $\mathcal{O}(md_i)$ for each $x$ and $i = 1,\ldots,k$.  Thus, with $n$ points making up the space, $\Omega = \bigcup_{i=1}^kV_i$, the complexity of constructing the collection $\{\xi_i\}_{i=1}^k$ is $\mathcal{O}(mnr)$, where $r = \sum_{i=1}^k d_i$.  

Similarly, one can expand $\zeta_i$ using the definition of $y_i$ to show that
$$\|y_i\|_2^2 \propto \left [ \|m_i^{(l)} - \bar{x}_i\|_2^2 + \left ( 2(1 - \alpha) + (1-\alpha)^2 \right )\|U_i^T(m_i^{(l)} - \bar{x}_i)\|_2^2 \right ],$$
\noindent which has complexity $\mathcal{O}(md_i)$.  Since $\|U_iU_i^Ty_i\|_2^2$ also has complexity $\mathcal{O}(md_i)$, we see that the cost of computing each $\zeta_i$ is $\mathcal{O}(mr)$, leading to an overall complexity of $\mathcal{O}(mnr)$ for the construction of the collection $\{\gamma_i\}_{i=1}^k$.  Since this matches the complexity of solving a problem instance using $\alpha \in \{0,1\}$, we see that the gradient descent steps can be performed with no additional cost in the case of $\alpha \in (0.1)$.

\subsection{Initialization for MOR}

Here, we briefly describe the initialization process used with our MOR experiments.  Let $\{x_i\}_{i=1}^n \subset \mathbb{R}^m$ be a collection  of model snapshots.  We assume each $x_i$ has mean zero.  We then compute the correlations, $\rho_{ij}$, between each point and take absolute values:
$$c_{ij} = |\rho_{ij}|.$$

\noindent Using $k$, the number of initial Voronoi sets prescribed by the user, we calculate the $k-1$ quantiles corresponding to the $c_{ij}$.  This will induce a partition of $[0,1]$ into $k$ disjoint sets.  After sorting the $c_{ij}$ in descending order, we associate each pair $x_i,x_j$ to the Voronoi set associated with the quantile partition they belong to.  Once a pair are assigned, they are removed from be considered again.  The process repeats until all snapshots are selected.  In the event that $n$ is odd, we assign the last point to the Voronoi set corresponding to the smallest quantile.


\end{document}